\documentclass[11pt]{amsart}

\usepackage{amsmath,amssymb,latexsym,amsthm}
\usepackage[all]{xy}

\usepackage{mathrsfs}

\newtheorem{theorem}{Theorem}[section]
\newtheorem{lemma}[theorem]{Lemma}
\newtheorem{corollary}[theorem]{Corollary}
\newtheorem{proposition}[theorem]{Proposition}

\theoremstyle{definition}
\newtheorem{definition}[theorem]{Definition}

\theoremstyle{remark}
\newtheorem{observation}[theorem]{Observation}
\newtheorem{remark}[theorem]{Remark}

\numberwithin{equation}{section}

\begin{document}

\title{Enriched Reedy categories}

\author{Vigleik Angeltveit}
\thanks{This research was partially conducted during the period the author was employed by the Clay Mathematics Institute as a Liftoff Fellow}
\address{University of Chicago \\ Department of Mathematics \\ 5734 S University Ave \\ Chicago IL 60608 \\ USA}
\email{vigleik@math.uchicago.edu}
\urladdr{http://www.math.uchicago.edu/$\sim$vigleik}

\subjclass[2000]{18G55}

\newcommand{\A}{\mathcal{A}}
\newcommand{\B}{\mathcal{B}}
\newcommand{\C}{\mathcal{C}}
\newcommand{\D}{\mathcal{D}}
\newcommand{\E}{\mathcal{E}}
\newcommand{\F}{\mathbb{F}}
\newcommand{\I}{\mathcal{I}}
\newcommand{\K}{\mathcal{K}}
\newcommand{\M}{\mathcal{M}}
\newcommand{\R}{\mathbb{R}}
\newcommand{\V}{\mathcal{V}}
\newcommand{\W}{\mathbb{W}}
\newcommand{\Z}{\mathbb{Z}}
\newcommand{\sma}{\wedge}
\newcommand{\lar}{\longrightarrow}
\newcommand{\sar}{\rightarrow}
\newcommand{\colim}{\varinjlim}
\newcommand{\hEn}{\widehat{E(n)}}
\newcommand{\bD}{\mathbf{\Delta}}
\newcommand{\noco}{\mathbf{\Delta} \Sigma}
\newcommand{\Top}{\mathcal{T}\!\mathit{op}}
\newcommand{\gup}{\overrightarrow{g}}
\newcommand{\gdown}{\overleftarrow{g}}

\def\arrow#1{\overset{#1}{\lar}}

\begin{abstract}
We define the notion of an enriched Reedy category, and show that if $\A$ is a $\C$-Reedy category for some symmetric monoidal model category $\C$ and $\M$ is a $\C$-model category, the category of $\C$-functors and $\C$-natural transformations from $\A$ to $\M$ is again a model category.
\end{abstract}

\maketitle

\section{Introduction}
A Reedy category is a category with a notion of an injective and a surjective morphism such that any morphism can be factored uniquely as a surjection followed by an injection. The simplicial indexing category is the prototypical example of a Reedy category, and if $\A$ is a Reedy category then so is the opposite category $\A^{op}$. A theorem of Dan Kan says that given a Reedy category $\A$ and a model category $\M$, the category $\M^\A$ of functors from $\A$ to $\M$ and natural transformations of such functors is again a model category, with the model structure described in Definition \ref{classical_reedy} below.

This should be compared to weak equivalences and fibrations in a diagram category as being defined levelwise, an approach that only works if $\M$ is cofibrantly generated, and weak equivalences and cofibrations being defined levelwise, which only works if $\M$ is combinatorial, a very strong condition to put on $\M$.

We are interested in an enriched version of this theory. Fix a symmetric monoidal model category $\C$. We will define a $\C$-Reedy category as a category which is enriched over $\C$ and satisfies a suitable analog of the unique factorization axiom, plus a cofibrancy condition. We prove that the category of $\C$-functors and $\C$-natural transformations from a $\C$-Reedy category $\A$ to a $\C$-model category $\M$ is a model category, and that something stronger is true: the functor category from $\A$ to $\M$, which is another category enriched over $\C$, is a $\C$-model category.

The results in this paper will be used to retain homotopical control in \cite{An_cyclic} and \cite{AnTHH}, where we define the cyclic bar construction on an $A_\infty$ $H$-space and use this to give a direct definition of topological Hochschild homology and cohomology of $A_\infty$ ring spectra in a way that is amenable to calculations.

This paper draws heavily on Hirschhorn's book \cite{Hi03}, particularly Chapter 15, and the author would like to thank Philip Hirschhorn for his help. The author would also like to thank Michael Shulman for reading an earlier version of the paper and finding several mistakes, and the referee for finding some more.

The author would also like to thank Justin Noel and Michael Shulman for noticing that the original definition of an enriched Reedy category (Definition \ref{d:enrichedReedy}) was imprecise. The updated version corrects this.

\section{Reedy categories}
We start by recalling a number of things from \cite{Hi03}.

\begin{definition} \label{def_reedy}
A Reedy category is a small category $\A$ together with two subcategories $\overrightarrow{\A}$ (the direct subcategory) and $\overleftarrow{\A}$ (the inverse subcategory), both of which contain all the objects of $\A$, together with a degree function assigning a nonnegative integer to each object in $\A$, such that
\begin{enumerate}
\item Every non-identity morphism of $\overrightarrow{\A}$ raises degree.
\item Every non-identity morphism of $\overleftarrow{\A}$ lowers degree.
\item Every morphism $g : \alpha \sar \beta$ in $\A$ has a unique factorization
\begin{equation}
\alpha \arrow{\gdown} \gamma \arrow{\gup} \beta
\end{equation}
with $\gdown$ a morphism in $\overleftarrow{\A}$ and $\gup$ a morphism in $\overrightarrow{\A}$.
\end{enumerate}
\end{definition}

The canonical example of a Reedy category is the cosimplicial indexing category $\bD$ with ordered sets $\mathbf{n}=\{0,1,\ldots,n\}$ and order-preserving maps. In this case $\overrightarrow{\bD}$ is the subcategory of injective maps and $\overleftarrow{\bD}$ is the subcategory of surjective maps.

Now let $\M$ be a model category, and suppose $X$ is a functor $\A \sar \M$. By model category we mean a closed model category, and we take as part of the definition that $\M$ is complete and cocomplete and that the factorizations into a cofibration followed by a trivial fibration, or a trivial cofibration followed by a fibration, are functorial. This is the version of Quillen's axioms found for example in \cite[Definition 7.1.3]{Hi03}.

\begin{definition}
Let $\alpha$ be an object in $\A$. The latching object $L_\alpha X$ is the colimit
\begin{equation}
L_\alpha X=\varinjlim_{\partial (\overrightarrow{\A}/\alpha)} X,
\end{equation}
where $\overrightarrow{\A}/\alpha$ is the category of objects over $\alpha$ and $\partial (\overrightarrow{\A}/\alpha)$ is the full subcategory containing all the objects except the identity on $\alpha$.

The matching object $M_\alpha X$ is the limit
\begin{equation}
M_\alpha X=\varprojlim_{\partial (\alpha/\overleftarrow{\A})} X,
\end{equation}
where $\alpha/\overleftarrow{\A}$ is the category of objects under $\alpha$ and $\partial (\alpha/\overleftarrow{\A})$ is the full subcategory containing all the objects except the identity on $\alpha$.
\end{definition}

\begin{remark} \label{factorLM}
An element in the direct limit system defining $L_\alpha X$ is a pair $(X_\beta, \beta \sar \alpha)$ and an element in the inverse limit system defining $M_\alpha X$ is a pair $(X_\gamma, \alpha \sar \gamma)$. Let $f$ be the composite $\beta \sar \alpha \sar \gamma$. Then we have a map $f_* :  X_\beta \sar X_\gamma$. It is not hard to check that this induces a map $L_\alpha X \sar M_\alpha X$, and that $X_\alpha$ provides a factorization of this map as $L_\alpha X \sar X_\alpha \sar M_\alpha X$.
\end{remark}

\begin{definition} \label{classical_reedy}
Let $X$ and $Y$ be functors $\A \sar \M$, and let $f : X \sar Y$ be a natural transformation.
\begin{enumerate}
\item The map $f$ is a Reedy weak equivalence if each 
\begin{equation}
f_\alpha : X_\alpha \lar Y_\alpha
\end{equation}
is a weak equivalence.
\item The map $f$ is a Reedy cofibration if each
\begin{equation}
X_\alpha \cup_{L_\alpha X} L_\alpha Y \lar Y_\alpha
\end{equation}
is a cofibration.
\item The map $f$ is a Reedy fibration if each
\begin{equation}
X_\alpha \lar Y_\alpha \times_{M_\alpha Y} M_\alpha X
\end{equation}
is a fibration.
\end{enumerate}
\end{definition}

We recall the following theorem, which is due to Dan Kan, from \cite[Theorem 15.3.4]{Hi03}:

\begin{theorem} \label{Reedy_model}
Let $\A$ be a Reedy category and let $\M$ be a model category. Then the category $\M^\A$ of functors from $\A$ to $\M$ with the Reedy weak equivalences, Reedy cofibrations and Reedy fibrations is a model category.

If $\M$ is a simplicial model category, then $\M^\A$ is again a simplicial model category.
\end{theorem}

\section{Enriched categories}
The purpose of this section is to introduce some notation and to recall some of the basic facts about enriched categories we will need. The canonical reference for enriched category theory is \cite{Kelly}.

Let $(\C,\otimes,I)$ be a closed symmetric monoidal category and let $\D$ be a category which is enriched over $\C$. Given objects $\alpha$ and $\beta$ in $\D$, we will write $Hom_\D(\alpha,\beta)$ for the $Hom$ object in $\C$ while using $hom_\D(\alpha,\beta)$ for the underlying $Hom$ set, defined by $hom_\D(\alpha,\beta)=hom_\C(I,Hom_\D(\alpha,\beta))$. We let $\D_0$ denote the underlying category of $\D$, so $Hom_{\D_0}(\alpha,\beta)=hom_\D(\alpha,\beta)$.

If $\D$ and $\E$ are enriched over $\C$, we write $hom(\D,\E)$ for the category of $\C$-functors and $\C$-natural transformations from $\D$ to $\E$. An object $X$ in $hom(\D,\E)$ consists of an object $X(\alpha)$ in $\E$ for each object $\alpha$ in $\D$, and a map $Hom_\D(\alpha,\beta) \sar Hom_\E(X(\alpha),X(\beta))$ in $\C_0$ for each pair $(\alpha,\beta)$ of objects in $\D$. Here $\C_0$ is the underlying category of $\C$, viewed as a category enriched over itself. A morphisms $F$ in $hom(\D,\E)$ from $X$ to $Y$ is a collection of maps $F_\alpha : I \sar Hom_\E(X(\alpha),Y(\alpha))$ in $\C$, or equivalently a collection of maps $X(\alpha) \sar Y(\alpha)$ in $\E_0$ satisfying certain compatibility conditions. This compatibility says that \cite[Diagram 1.7]{Kelly} is required to commute.

Sometimes it is also possible to define a category $Hom(\D,\E)$ which is enriched over $\C$ by imitating the description of $hom(\D,\E)(X,Y)$ for unenriched categories as an equalizer. The objects in $Hom(\D,\E)$ are the same as in $hom(\D,\E)$ but $Hom(\D,\E)(X,Y)$ is defined as the equalizer
\begin{multline}
Hom(\D,\E)(X,Y) \sar \prod_{\alpha \in \D} Hom_\E(X(\alpha),Y(\alpha)) \\ 
\rightrightarrows \prod_{\alpha,\beta \in \D} Hom_\C(Hom_\D(\alpha,\beta), Hom_\E(X(\alpha),Y(\beta)))
\end{multline}
if it exists. If $Hom(\D,\E)(X,Y)$ exists for all $\C$-functors $X$ and $Y$ from $\D$ to $\E$, then this defines a $\C$-category $Hom(\D,\E)$. This category is usually called the functor category from $\D$ to $\E$.

Now suppose that $\C$ is a monoidal model category, i.e., $\C$ is both a closed symmetric monoidal category and a model category, and these structures are compatible in the following sense: If $i:A \sar B$ and $j:K \sar L$ are cofibrations in $\C$, then the induced map
\begin{equation} \label{eq:pushoutproduct}
L \otimes A \cup_{K \otimes A} K \otimes B \lar L \otimes B
\end{equation}
is a cofibration that is trivial if either $i$ or $j$ is. This condition is called the pushout-product axiom.

We say that $\M$ is a $\C$-model category if $\M$ is a model category which is enriched, tensored and cotensored over $\C$ and satisfies the analog of the pushout-product condition, i.e., if $i:A \sar B$ is a cofibration in $\M$ and $j:K \sar L$ is a cofibration in $\C$, then the induced map (\ref{eq:pushoutproduct}) is a cofibration in $\M$ that is trivial if $i$ or $j$ is.

\begin{remark}
In \cite{Ho99} Hovey has another axiom which says that the canonical map $Q(I) \otimes X \sar I \otimes X \cong X$ is a weak equivalence for cofibrant $X$. This is important when passing to the homotopy category, but will not play a role here because we always work on the level of the model category.
\end{remark}

A monoidal model category $\C$ is sometimes called a Quillen ring, and a $\C$-model category is sometimes called a Quillen module.

If $\C$ is the category of simplicial sets then the pushout-product axiom is the extra condition that makes a model category which is enriched, tensored, and cotensored over $\C$ into a simplicial model category. If $\C$ is topological spaces, by which we mean compactly generated weak Hausdorff spaces, then a $\C$-category is a topological model category.

The pushout-product axiom has some immediate consequences. For example, it follows that if $A \sar B$ is a (trivial) cofibration in $\M$ then so is $K \otimes A \sar K \otimes B$ for any cofibrant $K$ in $\C$. Similarly, if $X \sar Y$ is a (trivial) fibration then so is $X^K \sar Y^K$ for any cofibrant $K$ in $\C$.

We will sometimes write $F(K,X)$ instead of $X^K$ for the cotensor.

\section{Enriched Reedy categories}
Next we define the notion of a $\C$-Reedy category. Here $\C$ is still a monoidal model category and $\M$ is a $\C$-model category. We bootstrap ourselves from the definition of a regular Reedy category.

\begin{definition} \label{d:enrichedReedy}
A $\C$-Reedy category is a small category $\A$ enriched over $\C$ together with an unenriched Reedy category $\B$ with the same objects and a decomposition
\begin{equation}
Hom_\A(\alpha,\beta)=\coprod_{g \in hom_\B(\alpha,\beta)} Hom_\A(\alpha,\beta)_g
\end{equation}
for each $Hom$ object such that the composition in $\A$ is induced by associative maps $Hom_\A(\beta,\gamma)_g \otimes Hom_\A(\alpha,\beta)_f \to Hom_\A(\alpha,\gamma)_{g \circ f}$ and the unit map $I \to Hom_\A(\alpha,\alpha)$ factors through an isomorphism to $Hom_\A(\alpha,\alpha)_{id}$. We require that the composition map
\begin{equation}
Hom_\A(\gamma,\beta)_{\gup} \otimes Hom_\A(\alpha,\gamma)_{\gdown} \lar Hom_\A(\alpha,\beta)_g
\end{equation}
is an isomorphism for each $g \in hom_\B(\alpha,\beta)$, where $g=\gup \circ \gdown$ is the factorization of $g$ in $\B$ as in Definition \ref{def_reedy}, together with the following cofibrancy condition. We set 
\begin{equation}
Hom_{\overrightarrow{\A}}(\alpha,\beta)=\coprod_{g \in Hom_{\overrightarrow{\B}}(\alpha,\beta)} Hom_\A(\alpha,\beta)_g
\end{equation}
and
\begin{equation}
Hom_{\overleftarrow{\A}}(\alpha,\beta)=\coprod_{g \in Hom_{\overleftarrow{\B}}(\alpha,\beta)} Hom_\A(\alpha,\beta)_g,
\end{equation}
and we require that each $Hom_{\overrightarrow{\A}}(\alpha,\beta)$ and $Hom_{\overleftarrow{\A}}(\alpha,\beta)$ is cofibrant in $\C$.
\end{definition}


Even though $\A$ has a discrete set of objects, the same is not true for $\partial(\overrightarrow{\A}/\alpha)$ and $\partial(\alpha/\overleftarrow{\A})$. Thus when defining the latching and matching object, we use the following enriched Kan extensions:

\begin{definition}
Let $X : \A \sar \M$ be a $\C$-functor. The latching object $L_\alpha X$ is the coequalizer
\begin{equation}
\coprod_{\beta<\gamma<\alpha} Hom_{\overrightarrow{\A}}(\gamma,\alpha) \otimes Hom_{\overrightarrow{\A}}(\beta,\gamma) \otimes X_\beta \rightrightarrows \coprod_{\beta<\alpha} Hom_{\overrightarrow{\A}}(\beta,\alpha) \otimes X_\beta \sar L_\alpha X,
\end{equation}
where one of the maps is given by the composition $Hom_{\overrightarrow{\A}}(\gamma,\alpha) \otimes Hom_{\overrightarrow{\A}}(\beta,\gamma) \sar \\ Hom_{\overrightarrow{\A}}(\beta,\alpha)$ and the other is given by $Hom_{\overrightarrow{\A}}(\beta,\gamma) \otimes X_\beta \sar X_\gamma$.

The matching object $M_\alpha X$ is the equalizer
\begin{multline}
M_\alpha X \sar \prod_{\beta<\alpha} F(Hom_{\overleftarrow{\A}}(\alpha,\beta),X_\beta) \\ 
\rightrightarrows \prod_{\beta<\gamma < \alpha} F(Hom_{\overleftarrow{\A}}(\gamma,\beta) \otimes Hom_{\overleftarrow{\A}}(\alpha,\gamma),X_\beta).
\end{multline}
\end{definition}

The category $\A$ has an obvious filtration, where $F^n \A$ is the full subcategory of $\A$ whose objects have degree less than or equal to $n$.

\begin{lemma}
(See Remark \ref{factorLM}.) Suppose $X$ is a functor $F^{n-1} \A$ $\sar \M$. Extending $X$ to a functor $F^n \A \sar \M$ is equivalent to choosing, for each object $\alpha$ of degree $n$, an object $X_\alpha$ and a factorization $L_\alpha X \sar X_\alpha \sar M_\alpha X$ of the natural map $L_\alpha X \sar M_\alpha X$.
\end{lemma}

\begin{proof}
This uses the unique factorization in the definition of a Reedy category, in the same way as in the proof of \cite[Theorem 15.2.1]{Hi03}.
\end{proof}

\begin{lemma}
Suppose that for every object $\beta$ of $\A$ of degree less than $\alpha$, the map $X_\beta \cup_{L_\beta X} L_\beta Y \sar Y_\beta$ is a (trivial) cofibration. Then $L_\alpha X \sar L_\alpha Y$ is a (trivial) cofibration.

Similarly, suppose that for every object $\beta$ of $\A$ of degree less than $\alpha$ the map $X_\beta \sar Y_\beta \times_{M_\beta Y} M_\beta X$ is a (trivial) fibration. Then $M_\alpha X \sar M_\alpha Y$ is a (trivial) fibration.
\end{lemma}

\begin{proof}
This is where we need the pushout-product axiom and that each $Hom_{\overrightarrow{\A}}(\alpha,\beta)$ and $Hom_{\overleftarrow{\A}}(\alpha,\beta)$ is cofibrant. We will do the case where each $X_\beta \cup_{L_\beta X} L_\beta Y \sar Y_\beta$ is a trivial cofibration, the other cases are similar. Let $E \sar B$ be a fibration. We have to show that any diagram
\begin{equation}
\xymatrix{ L_\alpha X \ar[r] \ar[d] & E \ar[d] \\ L_\alpha Y \ar@{.>}[ur] \ar[r] & B}
\end{equation}
has a lift. Classically we had to construct a map $Y_\beta \sar E$ for each object $\beta \sar \alpha$ in $\partial(\overrightarrow{\A}/\alpha)$ by induction on the degree of $\beta$. We need to make sure that these maps are compatible, so in our case we need to construct a map $Hom_{\overrightarrow{\A}}(\beta,\alpha) \otimes Y_\beta \sar E$ for each $\beta$ of degree less than $\alpha$.

We proceed by induction. Suppose we have maps $Hom_{\overrightarrow{\A}}(\gamma,\alpha) \otimes Y_\gamma \sar E$ for all $\gamma$ of degree less than $\beta$. We then have maps 
\begin{equation}
Hom_{\overrightarrow{\A}}(\beta,\alpha) \otimes Hom_{\overrightarrow{\A}}(\gamma,\beta) \otimes Y_\gamma \lar Hom_{\overrightarrow{\A}}(\gamma,\alpha) \otimes Y_\gamma \sar E
\end{equation}
for each $\gamma$ of degree less than $\beta$. These maps assemble to a map $Hom_{\overrightarrow{\A}}(\beta,\alpha) \otimes L_\beta Y \sar E$. We also have maps $Hom_{\overrightarrow{\A}}(\beta,\alpha) \otimes X_\beta \sar E$, so we get a diagram
\begin{equation}
\xymatrix{ Hom_{\overrightarrow{\A}}(\beta,\alpha) \otimes (X_\beta \cup_{L_\beta X} L_\beta Y) \ar[r] \ar[d] & E \ar[d] \\
Hom_{\overrightarrow{\A}}(\beta,\alpha) \otimes Y_\beta \ar[r] \ar@{.>}[ur] & B}
\end{equation}
By assumption, each map $X_\beta \cup_{L_\beta X} L_\beta Y \sar Y_\beta$ is a trivial cofibration, and by the pushout-product axiom this remains true after tensoring with $Hom_{\overrightarrow{\A}}(\beta,\alpha)$, so we have a lift. These lifts are clearly compatible, and induce a lift $L_\alpha Y \sar E$.
\end{proof}

\begin{lemma}
A map $X \sar Y$ is a trivial Reedy cofibration if and only if each $X_\alpha \cup_{L_\alpha X} L_\alpha Y \sar Y_\alpha$ is a trivial cofibration.

Similarly, $X \sar Y$ is a trivial Reedy fibration if and only if each $X_\alpha \sar Y_\alpha \times_{M_\alpha Y} M_\alpha X$ is a trivial fibration.
\end{lemma}

\begin{proof}
We will only do the first part, the second part is dual. Recall that the pushout of a trivial cofibration is a trivial cofibration. Suppose that $f : X \sar Y$ is a trivial Reedy cofibration. We need to prove that each $X_\alpha \cup_{L_\alpha X} L_\alpha Y \sar Y_\alpha$ is a weak equivalence. By induction we can assume that $X_\beta \cup_{L_\beta X} L_\beta Y \sar Y_\beta$ is a weak equivalence for $\beta<\alpha$. By the previous lemma $L_\alpha X \sar L_\alpha Y$ is a trivial cofibration, so when we take the pushout over the map $L_\alpha X \sar X_\alpha$ we find that the map $X_\alpha \sar X_\alpha \cup_{L_\alpha X} L_\alpha Y$ is a trivial cofibration. By assumption the composite $X_\alpha \sar X_\alpha \cup_{L_\alpha X} L_\alpha Y \sar Y_\alpha$ is a weak equivalence, so by the two out of three axiom so is $X_\alpha \cup_{L_\alpha X} L_\alpha Y \sar Y_\alpha$.

The converse is similar.
\end{proof}

\begin{theorem} \label{Reedy_still_works}
Let $\A$ be a $\C$-Reedy category and let $\M$ be a $\C$-model category. Then the category $hom(\A,\M)$ of $\C$-functors and $\C$-natural transformations from $\A$ to $\M$ with the Reedy weak equivalences, Reedy cofibrations and Reedy fibrations is a model category.
\end{theorem}

\begin{proof}
If we have a diagram
\begin{equation}
\xymatrix{ A \ar[r] \ar[d] & X \ar[d] \\ B \ar[r] \ar@{.>}[ur] & Y}
\end{equation}
where $i : A \sar B$ is a Reedy cofibration and $p: X \sar Y$ is a Reedy fibration, with either $i$ or $p$ a weak equivalence, we need to construct a lift. We can do this by induction on the degree, using the diagrams
\begin{equation}
\xymatrix{ A_\alpha \cup_{L_\alpha A} L_\alpha B \ar[r] \ar[d] & X_\alpha \ar[d] \\
B_\alpha \ar[r] \ar@{.>}[ur] & Y_\alpha \times_{M_\alpha Y} M_\alpha X}
\end{equation}
and the previous lemma.
\end{proof}

\begin{theorem}
Let $\A$ be a $\C$-Reedy category and let $\M$ be a $\C$-model category. Then the functor category $Hom(\A,\M)$ is a $\C$-model category.
\end{theorem}

\begin{proof}
First of all, the category $Hom(\A,\M)$ exists because $\A$ is small and one of the axioms for a model category is that it has all small limits. In particular, the equalizer defining $Hom(\A,\M)(X,Y)$ for each $X$ and $Y$ exists.

We define $K \otimes X$ and $X^K$ for a $\C$-functor $X : \A \lar \M$ and an object $K \in \C$ objectwise, and it is clear that
\begin{equation}
Hom(K \otimes X,Y) \cong Hom(K,Hom(X,Y)) \cong Hom(X,Y^K)
\end{equation}
because this holds objectwise.

It remains to show that the pushout-product axiom holds. If $i : A \sar B$ is a Reedy cofibration in $hom(\A,\M)$ and $j : K \sar L$ is a cofibration in $\C$, we need to show that
\begin{equation}
L \otimes A \coprod_{K \otimes A} K \otimes B \lar L \otimes B
\end{equation}
is a Reedy cofibration in $hom(\A,\M)$. But this is equivalent to each
\begin{equation}
(L \otimes A \coprod_{K \otimes A} K \otimes B)_\alpha \coprod_{L_\alpha(L \otimes A \coprod_{K \otimes A} K \otimes B)} L_\alpha(L \otimes B) \lar (L \otimes B)_\alpha
\end{equation}
being a cofibration. By using that colimits commute with tensors, this is equivalent to each
\begin{equation}
(K \otimes B_\alpha) \coprod_{K \otimes (L_\alpha B \coprod_{L_\alpha A} A_\alpha)}L \otimes (L_\alpha B \coprod_{L_\alpha A} A_\alpha) \lar L \otimes B_\alpha
\end{equation}
being a cofibration, and this follows from the pushout-product axiom for $\M$. The case where $i$ or $j$ is also a weak equivalence is similar.
\end{proof}

\section{Homotopy limits and colimits}
Given $\C$-functors $A:\A \sar \M$ and $K:\A^{op} \sar \C$ we define $K \otimes_\A A$ as the coequalizer
\begin{equation}
\coprod_{\alpha,\beta \in \A} K_\beta \otimes Hom_\A(\alpha,\beta) \otimes A_\alpha \rightrightarrows \coprod_{\alpha \in \A} K_\alpha \otimes A_\alpha \sar K \otimes_{\A} A.
\end{equation}
Similarly, if $K:\A \sar \C$ we define $hom_\A(K,A)$ as the equalizer
\begin{equation}
hom_\A(K,A) \sar \prod_{\alpha \in \A} F(K_\alpha,A_\alpha) \rightrightarrows \prod_{\alpha,\beta \in \A} F(K_\alpha \otimes Hom_\A(\alpha,\beta),A_\beta).
\end{equation}

If $K$ is Reedy cofibrant we think of $K \otimes_\A A$ as a model for the homotopy colimit of $A$ and $hom_\A(K,A)$ as a model for the homotopy limit. In particular, if $\A$ is an enriched version of the simplicial indexing category then this gives a good notion of geometric realization (for suitable $K$).

\begin{theorem} (Compare \cite[Theorem 18.4.11]{Hi03}.) \label{thm:pushoutproduct}
Let $\A$ be a $\C$-Reedy category and let $\M$ be a $\C$-model category. If $j : A \sar B$ is a Reedy cofibration in $hom(\A,\M)$ and $i : K \sar L$ is a Reedy cofibration in $hom(\A^{op},\C)$, then
\begin{equation}
L \otimes_\A A \coprod_{K \otimes_\A A} K \otimes_\A B \lar L \otimes_\A B
\end{equation}
is a cofibration in $\M$ that is a weak equivalence if either $i$ or $j$ is.

Dually, if $p: X \sar Y$ is a Reedy fibration in $hom(\A,\M)$ and $i : K \sar L$ is a Reedy cofibration in $hom(\A,\C)$ then
\begin{equation}
hom_\A(L,X) \lar hom_\A(K,X) \times_{hom_\A(K,Y)} hom_\A(L,Y)
\end{equation}
is a fibration in $\M$ that is a weak equivalence if either $i$ or $p$ is.
\end{theorem}

\begin{corollary}
If $K$ is a Reedy cofibrant object in $hom(\A^{op},\C)$ and $f:X \sar Y$ is a weak equivalence of Reedy cofibrant objects in $hom(\A,\M)$, then the induced map $f_* : K \otimes_\A X \sar K \otimes_\A Y$ is a weak equivalence of cofibrant objects in $\M$.

Dually, if $K$ is a Reedy cofibrant object in $hom(\A,\C)$ and $f:X \sar Y$ is a weak equivalence of Reedy fibrant objects in $hom(\A,\M)$ then the induced map $f^* : hom_\A(K,X) \sar hom_\A(K,Y)$ is a weak equivalence of fibrant objects in $\M$.
\end{corollary}

\section{The Reedy category $\A_P$}
Let $\noco$ be the category of noncommutative sets, as in \cite{PiRi}. The objects in $\noco$ are finite sets $\mathbf{n}=\{0,1,\ldots,n\}$ and the morphisms are maps of finite sets together with a linear ordering of each inverse image of an element. Now let $\A$ be a Reedy category over $\noco$, i.e., $\A$ comes with a functor $U : \A \sar \noco$. Also let $P$ be an operad, by which we mean non-$\Sigma$ operad, in $\C$ with $P(0)=P(1)=I$ and each $P(n)$ cofibrant.

As in \cite[Definition 3.1]{An_cyclic}, we define a new category $\A_P$ enriched over $\C$ as follows. The objects are the same as in $\A$, but the $Hom$ objects are given by
\begin{equation}
Hom_{\A_P}(\alpha,\beta)=\coprod_{f \in Hom_\A(\alpha,\beta)} P[f],
\end{equation}
where $P[f]= \bigotimes_{i \in U \beta} P(Uf^{-1}(i))$. Composition in $\A_P$ is defined using the structure maps for $P$.

\begin{proposition} \label{prop:APisReedy}
The category $\A_P$ is a $\C$-Reedy category.
\end{proposition}

\begin{proof}
The decomposition of $Hom_{\A_P}(\alpha,\beta)$ as a coproduct over $Hom_\A(\alpha,\beta)$ is the one in the above definition. The condition $P(0)=P(1)=I$ ensures that the direct subcategory $\overrightarrow{\A}_P$ is in fact equal to $\overrightarrow{\A}$, and it is easy to see that $\A_P$ satisfies the unique factorization condition. The condition that each $P(n)$ is cofibrant ensures that the cofibrancy hypothesis in the definition is satisfied.
\end{proof}

\begin{corollary}
Let $\M$ be a $\C$-model category. Then the category $hom(\A_P,\M)$ of $\C$-functors and $\C$-natural transformations from $\A_P$ to $\M$ is a model category, and the functor category $Hom(\A_P,\M)$ is a $\C$-model category.

Similarly, $hom(\A_P^{op},\M)$ is a model category and $Hom(\A_P^{op},\M)$ is a $\C$-model category.
\end{corollary}

From the proof of Proposition \ref{prop:APisReedy} we observe the following:

\begin{observation} \label{samelatching}
Because the direct subcategory does not change when we pass from $\A$ to $\A_P$, if $\A$ is isomorphic to $\bD^{op}$ and $X : \A_P \sar \M$ the usual description of latching objects as the coequalizer
\begin{equation}
\coprod_{0 \leq i<j \leq n-1} X_{n-2} \rightrightarrows \coprod_{0 \leq i \leq n-1} X_{n-1} \lar L_n X
\end{equation}
as in \cite[Proposition 15.2.6]{Hi03} is still valid.

Dually, the usual description of matching objects for $Y : \A_P^{op} \sar \M$ does not change when we pass from $\A$ to $\A_P$.
\end{observation}

We have two examples of Reedy categories over $\noco$ which are isomorphic to $\bD^{op}$. Let $^{01} \! \bD$ be (a skeleton of) the category of doubly based totally ordered sets. The objects are totally ordered sets of cardinality at least $2$, and the morphisms are order-preserving maps which preserve the minimal and maximal element. For the second example, let $^0 \! \bD C$ be the category whose objects are cyclically ordered sets with a given basepoint, and whose morphisms are maps of cyclically ordered sets which preserve the basepoint.

\begin{lemma} (\cite[Lemma 3.3]{An_cyclic}.)
The categories $^{01} \! \bD$ and $^0 \! \bD C$ are isomorphic to $\bD^{op}$.
\end{lemma}

\section{The associahedra operad}
Now let $\C$ be either simplicial sets or topological spaces, and let $P=\K$ be the associahedra operad in $\C$, see \cite{St63} and \cite{An_cyclic}. Also let $\A$ be either $^{01} \! \bD$ or $^0 \! \bD C$, so $\A \cong \bD^{op}$. Then $\A_\K$ is a $\C$-Reedy category. Also recall from \cite{An_cyclic} the definition of geometric realization $|X|$ for a $\C$-functor $\A_\K \sar \M$ as $\K \otimes_{\A_\K} X$ if $\A={}^{01} \! \bD$ and $\mathcal{W} \otimes_{\A_\K} X$ if $\A={}^0 \! \bD C$, where $\mathcal{W}$ is the cyclohedra.

\begin{proposition} \label{KReedycofibrant}
The $\C$-functor $\K : {}^{01} \! \bD_\K^{op} \sar \C$ is Reedy cofibrant. Similarly, $\mathcal{W} : {}^0 \! \bD C_\K \sar \C$ is Reedy cofibrant.
\end{proposition}

\begin{proof}
Let the degree function for the Reedy category $^{01} \! \bD_\K$ be the one sending a set with $n+2$ elements to $n$, so it corresponds to the standard degree function on $\bD^{op}$ under the isomorphism $^{01} \! \bD \cong \bD^{op}$.

Then we need to check that each $L_n \K \sar K_{n+2}$ is a cofibration. But $L_n \K \cong \partial K_{n+2}$, the union of the faces of $\K$, so $L_n \K \sar K_{n+2}$ is homeomorphic to $S^{n-1} \sar D^n$, which is certainly a cofibration.

The other case is similar.
\end{proof}

\begin{corollary}
If $f:X \sar Y$ is a Reedy weak equivalence between Reedy cofibrant $\C$-functors $\A_\K \sar \M$ with $\A={}^{01} \! \bD$ or $^0 \! \bD C$, then $f$ induces a weak equivalence $f_* :  |X| \sar |Y|$.
\end{corollary}

\begin{proof}
This follows from Theorem \ref{thm:pushoutproduct} and Proposition \ref{KReedycofibrant}
\end{proof}

We also get the expected spectral sequences in this setup.

\begin{theorem} \label{skeletalfilt}
Let $\M$ be a pointed $\C$-model category, let $X: {}^{01} \! \bD_\K \sar \M$ or $^0 \! \bD C_\K \sar \M$ be Reedy cofibrant and let $E$ be a homology theory. Then the skeletal filtration gives a spectral sequence
\begin{equation}
E^2_{p,q}=H_p(E_q(X)) \Longrightarrow E_{p+q} |X|.
\end{equation}
\end{theorem}

\begin{proof}
The proof is similar to the classical case. By Proposition \ref{KReedycofibrant} and our definition of geometric realization it follows that each $sk_{n-1} X \sar sk_n X$ is a cofibration in $\M$. To build the spectral sequence we only have to identify the filtration quotients and the $d^1$-differential.

Each filtration quotient looks like $K_{n+2}/\partial K_{n+2} \otimes X_n/L_n X$ in the first case and $W_{n+1}/\partial W_{n+1} \otimes X_n/L_n X$ in the second case, and this identifies the $E^1$-term as the normalized chain complex associated to the graded simplicial abelian group $E_* X$. The identification of the $E^2$-term is standard.
\end{proof}

There is also a dual setup for Reedy fibrant right modules.

\begin{theorem} \label{totfilt}
Let $Y$ be a Reedy fibrant functor $^{01} \! \bD_\K \sar \M$ or $^0 \! \bD C_\K \sar \M$, and let $E$ be a homology theory. Then the total space filtration gives a spectral sequence
\begin{equation}
E_2^{p,q}=H^p(E_q(Y)) \Longrightarrow E_{q-p} Tot(Y).
\end{equation}
\end{theorem}

While the spectral sequence coming from the skeletal filtration usually has good convergence properties, we need additional conditions to guarantee convergence of the spectral sequence coming from the total object filtration. See for example \cite{Bo89} for details.

\bibliographystyle{amsplain}

\bibliography{b2}

\end{document}